\documentclass{amsart}
\usepackage{amscd}
\usepackage{amsfonts}
\usepackage{amsmath}
\usepackage{amssymb}
\usepackage{latexsym}
\usepackage{amsthm,color}
\usepackage[all]{xy}
\usepackage{epsfig}
\usepackage{graphicx}
\usepackage{color}
\newtheorem{lemma}{Lemma}[section]

\newtheorem{theorem}{Theorem}
\newtheorem{corollary}{Corollary}

\theoremstyle{definition}
\newtheorem{definition}{Definition}
\theoremstyle{definition}
\newtheorem{example}{Example}[section]
\theoremstyle{definition}
\newtheorem{remark}{Remark}[section]

\newcommand{\C}{\mathbb{C}}

\newcommand{\R}{\mathbb{R}}

\begin{document}
\title[
%The concept of duality in and 
Characterization of Legendre involution  
%on the class of hyperplane families  creating envelopes %frontals
]
{A characterization of the Legendre involution  \\ 
on the class of generic frontals 
%envelopes of 
%hyperplane families % creating envelopes 
}
%The concept of duality in frontal geometry, \\ 
%and the characterization of the Legendre involution 
%}
%%%%%%%%%%%%%%%%%%%%%%%%%%%%%%%%%%%%%%%%%%%%%%%%%%%%%
%%%%%%%%%%%%%%%%%%%%%%%%%%%%%%%%%%%%%%%%%%%%%%  
\author[T.~Nishimura]{Takashi Nishimura%\thanks{Corresponding author.} 
}
\address{
Research Institute of Environment and Information Sciences,  
Yokohama National University, 
240-8501 Yokohama, JAPAN}
\email{nishimura-takashi-yx@ynu.ac.jp}
%%%%%%%%%%%%%%%%%%%%%%%%%%%%%%%%%%%%%%%%%%%%%%%  
%%%%%%%%%%%%%%%%%%%%%%%%%%%%%%%%%%%%%%%%%%%%%%% 
\begin{abstract}
%In the main theorem of this paper, 
We show that,  
under an additional mild assumption, 
on the class of generic frontals,  
%envelopes of hyperplane families with generic 
%Gauss mapping, 
any involution 
whose fixed point set is exactly the same as 
the fixed point set of the Legendre involution 
%\[
%\mathcal{GY}=\left\{
%\left(\left.
%\left(\theta_1, \ldots, \theta_n\right), 
%\frac{1}{2}\left(\theta_1^2+\cdots +\theta_n^2\right), 
%\left(\theta_1, \ldots, \theta_n\right)
%\right): \mathbb{R}^n\to \mathbb{R}^{2n+1} 
%\; \right|\; 
%\theta_i: \mathbb{R}^n\to \mathbb{R}\quad C^\infty
%\right\}, 
%\] 
must be the Legendre involution (Theorem \ref{theorem2} in \S 1).     
Moreover, its natural complexification (Corollary \ref{corollary1} in \S 3) 
is simultaneously shown.   
%; 
%where $\overline{ 
%\mbox{reg}\left(\theta_1, \ldots, \theta_n\right)}$ 
%stands for the topological closure 
%of the set consisting of regular points of %the $C^\infty$ mapping 
%$\left(\theta_1, \ldots, \theta_n\right): \mathbb{R}^n\to \mathbb{R}^n$.
\end{abstract}
 %, 53A40, 53A04}

%\thanks{}

\date{}

\maketitle

\section{Introduction\label{section1}}
Throughout this paper, let $n$ be a positive integer.    
Moreover, all functions and mappings are of class $C^\infty$ 
unless otherwise stated.   
\par 
\medskip 
The purpose of this paper is, on the class of generic frontals, 
obtaining a corresponding version of the following 
breakthrough result due to Shiri Artstein-Avidan and Vitali Milman.   
\begin{theorem}[\cite{artsteinavidanmilman}] \label{theorem1}
Assume a transform 
$\widetilde{\mathcal{T}} : Cvx \left(\mathbb{R}^n\right)\to 
Cvx \left(\mathbb{R}^n\right)$ 
satisfies 
\begin{enumerate}
\item $\widetilde{\mathcal{T}}\circ\widetilde{\mathcal{T}}(\phi) = \phi$ 
for any $\phi\in Cvx \left(\mathbb{R}^n\right)$,  
\item $\phi\le \psi$   implies 
$\widetilde{\mathcal{T}}(\phi) \ge 
\widetilde{\mathcal{T}}(\psi)$ for any 
$\phi, \psi\in Cvx \left(\mathbb{R}^n\right)$.   
\end{enumerate}
Then, $\widetilde{\mathcal{T}}$ 
is essentially the classical Legendre transform 
$\widetilde{\mathcal{L}}: Cvx\left(\mathbb{R}^n\right)\to 
Cvx\left(\mathbb{R}^n\right)$; 
namely there exists
a constant $C_0\in \mathbb{R}$, a vector $v_0\in \mathbb{R}^n$, 
and an invertible symmetric linear
transformation $B\in GL_n$ such that
\[
(\widetilde{\mathcal{T}}(\phi))(x) 
= (\widetilde{\mathcal{L}}(\phi))(Bx + v_0) + 
\left\langle x, v_0\right\rangle + C_0.   
\]
\end{theorem}
%%%%%%%%%%%%%%%%%%%%%%%%%%%%%%%%%%%%%%%%%%%  
\noindent 
Here, $Cvx\left(\mathbb{R}^n\right)$, 
$\langle{\;}, {\;}\rangle$ and 
$\widetilde{\mathcal{L}}: Cvx\left(\mathbb{R}^n\right)\to 
Cvx\left(\mathbb{R}^n\right)$ 
stand for the set consisting of lower semi-continuous 
convex functions  
$\phi: \mathbb{R}^n\to \mathbb{R}\cup \left\{\pm \infty\right\}$, 
the standard scalar product on $\mathbb{R}^n$ and 
the classical Legendre transform defined by 
$
\left(\widetilde{\mathcal{L}}(\phi)\right)(x)=
\sup_{y\in \mathbb{R}^n}
\left(\langle x, y\rangle - \phi(y)\right)
$ 
respectively.   
The Legendre transform $\widetilde{\mathcal{L}}
: Cvx\left(\mathbb{R}^n\right)\to 
Cvx\left(\mathbb{R}^n\right)$ is a very important notion 
having many applications 
(see for example, \cite{arnoldmechanics, hormander, rockafellar}).    
%%%%%%%%%%%%%%%%%%%%%%%%%%%%%%%%%%%%%%%%%%%%%%%% 
\par 
\smallskip 
%We want to obtain a generic frontal version of Theorem \ref{theorem1}.   
A frontal, which is a generalized notion of wavefront, 
was started to study relatively recently.   
Thus, although wavefronts seem to be already well-studied 
(see for instance, \cite{arnold, arnoldetal}), frontals have been  
less studied so far.    As for frontals, there is a well-written 
expository article \cite{ishikawa} which is recommended to readers.    
A {\it frontal} 
is a $C^\infty$ mappping $\varphi: \mathbb{R}^n\to \mathbb{R}^{n+1}$ 
such that there exists a $C^\infty$ mapping 
$\nu: \mathbb{R}^n\to S^n$ 
(called the \textit{Gauss mapping of} $\varphi$)
satisfying $\langle d\varphi_x(v), \nu(x)\rangle=0$ 
for any $x\in \mathbb{R}^n$ and $v\in T_x\mathbb{R}^n$ where $S^n$ is the 
uniit sphere in $\R^{n+1}$;  
%and no restriction on 
%the Legendre mapping 
%$(g, \nu): \mathbb{R}^n\to \mathbb{R}^{n+1}\times S^n$.    
while a {\it wavefront} is a $C^\infty$ mappping 
${\varphi}: \mathbb{R}^n\to \mathbb{R}^{n+1}$ 
such that there exists a $C^\infty$ mapping 
${\nu}: \mathbb{R}^n\to S^n$ satisfying 
$\langle d{\varphi}_x(v), {\nu}(x)\rangle=0$ 
for any $x\in \mathbb{R}^n$ and $v\in T_x\mathbb{R}^n$ and 
the mapping $({\varphi}, {\nu}): 
\mathbb{R}^n\to \mathbb{R}^{n+1}\times S^n$ 
(called the {\it Legendre mapping}) 
is immersive.    Thus, the difference between a wavefront and 
a frontal is only one point whether the Legendre mapping is 
non-singular or not necessarily no-singular.   
For instance, consider the mapping $\varphi: \mathbb{R}\to \mathbb{R}^2$ 
%of one variable 
%$f_1(x)=x^5$ and its graph 
defined by $\varphi(x)=\left(x^2, x^5\right)$.    Then, it is easily 
seen that $\varphi$ is a frontal but not a wavefront.    
Since the class on which the Legendre involution is considered 
should be wide as much as possible, in this paper we deal with 
frontals rather than wavefronts.   
%do not restrict ourselves to wavefronts,   
By using the results on envelopes created by 
hyperplane families obtained in \cite{nishimura} 
(see also \cite{nishimurahmj} which is, concentrated on the case $n=1$, an 
easy to understand expository article on the envelopes created by 
line families in the plane), 
it is possible to construct a natural bijective mapping 
$\mathcal{C}$ from the set $\mathcal{F}$ consisting of  
frontals $\varphi: \mathbb{R}^n\to \mathbb{R}^{n+1}$ 
%(For a frontal, \lq\lq being generic\rq\rq 
%is defined in Definition \ref{definition1}, which appears immediately
%after) 
to the set $\mathcal{X}$ consisting of $C^\infty$ mappings  
$\Phi: \mathbb{R}^n\to \mathbb{R}^{2n+1}$ given by 
\begin{eqnarray*} 
\Phi(x) & = & \left(\theta(x), a(x), b(x)\right), \\ 
\theta(x) & = & \left(\theta_1(x), \ldots, \theta_n(x)\right),  \\
b(x) & = & \left(b_1(x), \ldots, b_n(x)\right) 
\end{eqnarray*}
satisfying %the following three conditions:.   
%\begin{enumerate}
%\item 
%The following equality of one-forms $(*)$ is satisfied.   
\[
da = \sum_{i=1}^n b_id \theta_i,    
\leqno{(*)}
\]
where $da$ stands for the exterior derivative of the function 
$a: \mathbb{R}^n\to \mathbb{R}$.   
The condition $(*)$ is called the \textit{creative condition}.  
For details on the bijective mapping 
$\mathcal{C}: \mathcal{F}\to \mathcal{X}$, see Subsection \ref{subsection2.2}.     
%\item The set consisting of regular points of 
%$\left(\theta_1(x), \ldots, \theta_n(x)\right): \R^n\to \R^n$ 
%(denoted by $Reg(\theta)$) is dense in $\R^n$.   
%\item The set consisting of regular points of 
%$\left(b_1(x), \ldots, b_n(x)\right): \R^n\to \R^n$ 
%(denoted by $Reg(b)$) is dense in $\R^n$.   
%\end{enumerate}
\begin{definition}\label{definition1}
Given a frontal $\varphi\in \mathcal{F}$, 
set $\mathcal{C}(\varphi)=\left(\left(\theta_1, \ldots, \theta_n\right), 
a, \left(b_1, \ldots, b_n\right)\right)$, where 
$\mathcal{C}: \mathcal{F}\to \mathcal{X}$ 
is the above bijective 
mapping.    Then, $\varphi$ is said to be 
\textit{generic} if both of the set of regular points of the mapping 
$\left(\theta_1, \ldots, \theta_n\right): \mathbb{R}^n\to \mathbb{R}^n$  
 (denoted by $Reg(\theta)$) and the set of regular points of the mapping 
$\left(b_1, \ldots, b_n\right): \mathbb{R}^n\to \mathbb{R}^n$ 
(denoted by $Reg(b)$) 
are dense in $\R^{n}$.   
\end{definition}
\noindent
%Then, there is a natural bijective mapping 
%one-to-one correspondence between 
Let $\mathcal{GF}$ be the set consisting of generic frontals 
$\varphi: \R^n\to \R^{n+1}$ and let $\mathcal{GX}$ be the 
following set.      
\[ 
\mathcal{GX}  =   
\left\{\Phi=\left(\left(\theta_1, \ldots, \theta_n\right), a, 
\left(b_1, \ldots, b_n\right)\right): \mathbb{R}^n\to \mathbb{R}^{2n+1} 
\; \left|\; 
\begin{array}{c}
\theta_i, a, b_i: \mathbb{R}^n\to \mathbb{R}\;\; C^\infty\;\; 
(\forall i\in \{1, \ldots, n\}), \\ 
(*) \mbox{ is satisfied},  \\ 
Reg(\theta), Reg(b) \mbox{ are dense in $\mathbb{R}^n$}
\end{array}
\right.\right\}.    
\]
%For details on this correspondence, see Subsection \ref{subsection2.2}.   
%Under this correspondence, any transform of the set of 
%generic frontals %$\mathcal{F}\to \mathcal{F}$ 
%is uniquely identified with the corresponding transform 
%$\mathcal{GX}\to \mathcal{GX}$ and vice versa.   
By definition, it is clear that 
the restriction of $\mathcal{C}$ to $\mathcal{GF}$ is 
a bijective mapping 
$\mathcal{C}|_{\mathcal{GF}}: \mathcal{GF}\to \mathcal{GX}$.     
\begin{definition}[Appendix 4 in \cite{arnoldmechanics}]
{\rm 
The mapping $\mathcal{L}: \mathcal{X}\to 
C^\infty(\mathbb{R}^n, \mathbb{R}^{2n+1})$ defined by 
\[
\mathcal{L}\left(\left(\theta_1, \ldots, \theta_n\right),\; a,\;  
\left(b_1, \ldots, b_n\right)\right) 
= 
\left(\left(b_1, \ldots, b_n\right),\; \sum_{i=1}^n b_i\theta_i -a,\;  
\left(\theta_1, \ldots, \theta_n\right)\right)
\] 
is called the \textit{Legendre involution}.   
}
\end{definition}
\noindent 
Notice that in \cite{arnoldmechanics} 
the Legendre involution $\mathcal{L}$ is defined               
over a contact space.    On the other hand, 
in this paper, thanks to the creative condition $(*)$, 
no contact spaces are required.    
This viewpoint is meaningful (for details on this viewpoint, see     
%For details on this notice, see Remark \ref{remark2} in Subsection 
Subsection \ref{subsection2.1}.   
It is easily seen that the Legendre involution 
$\mathcal{L}$ is actually an involution $\mathcal{X}\to \mathcal{X}$.   
Moreover, even if $\mathcal{X}$ is replaced with $\mathcal{GX}$, again 
we have the property 
$\mathcal{L}\left(\mathcal{GX}\right)=\mathcal{GX}$.       
We want to obtain the characterization of the Legendre involution 
$\mathcal{L}: \mathcal{GX}\to \mathcal{GX}$ which corresponds to Theorem 
\ref{theorem1}.       
In order to obtain the condition corresponding to the order-reversing condition 
(2) of Theorem \ref{theorem1}, 
we pay attention to the fixed point set of the Legendre involution 
$\mathcal{L}: \mathcal{GX}\to \mathcal{GX}$.      
It is easily seen that the fixed point set of the Legendre involution 
$\mathcal{L}: \mathcal{X}\to \mathcal{X}$ and 
$\mathcal{L}: \mathcal{GX}\to \mathcal{GX}$ 
is exactly 
the following set $\mathcal{Y}$ and $\mathcal{GY}$ respectively.     
 \begin{eqnarray*}
\mathcal{Y} & = & \left\{
\left(
\left(\theta_1, \ldots, \theta_n\right), 
\frac{1}{2}\left(\theta_1^2+\cdots +\theta_n^2\right), 
\left(\theta_1, \ldots, \theta_n\right)
\right): \mathbb{R}^n\to \mathbb{R}^{2n+1}\; 
C^\infty
\right\}.   \\ 
%\; \right|\; 
%\theta_i: \mathbb{R}^n\to \mathbb{R}\; C^\infty
%, \; 
%\overline{ 
%Reg\left(\theta\right) \mbox{ is dense in }\mathbb{R}^n 
%= \mathbb{R}^n  
%\right\}.   \\ 
\mathcal{GY} & = & \left\{
\left(\left.
\left(\theta_1, \ldots, \theta_n\right), 
\frac{1}{2}\left(\theta_1^2+\cdots +\theta_n^2\right), 
\left(\theta_1, \ldots, \theta_n\right)
\right): \mathbb{R}^n\to \mathbb{R}^{2n+1}\; 
C^\infty
\; \right|\; 
%\theta_i: \mathbb{R}^n\to \mathbb{R}\; C^\infty
%, \; 
%\overline{ 
Reg\left(\theta\right) \mbox{ is dense in }\mathbb{R}^n 
%= \mathbb{R}^n  
\right\}.   
\end{eqnarray*}
Notice that $\mathcal{Y}$ (resp., $\mathcal{GY}$) 
is a proper subset of $\mathcal{X}$ (resp., $\mathcal{GX}$).    
%For more details on the Legendre involution 
%$\mathcal{L}: \mathcal{GX}\to \mathcal{GX}$ as well, see \S 2.  
\par 
The main result of this paper is 
as follows.   
%By now, we can state  
\begin{theorem}\label{theorem2} 
Assume a transform 
$\mathcal{T} : \mathcal{GX}\to 
\mathcal{GX}$ 
satisfies 
\begin{enumerate}
\item $\mathcal{T}\circ\mathcal{T}(\Phi) = \Phi$ 
for any $\Phi\in \mathcal{GX}$,  
\item For any $\Phi\in \mathcal{GX}$, 
$\mathcal{T}(\Phi)=\Phi$  if $\Phi\in \mathcal{GY}$ and 
$\mathcal{T}(\Phi) \ne \Phi$ if $\Phi\not\in \mathcal{GY}$.     
%$\mathcal{GY}$.   
\end{enumerate}
Assume moreover that there exists an analytic mapping 
$F: \mathbb{R}^{2n+1}\to \mathbb{R}^{2n+1}$ satisfying 
$\mathcal{T}(\Phi)=F\circ \Phi$ for any $\Phi\in \mathcal{GX}$.    
Then, $\mathcal{T}$ must be the Legendre involution 
$\mathcal{L}: \mathcal{GX}\to \mathcal{GX}$.     
\end{theorem}   
In \S 3, it turns out that there exists a transform 
$\mathcal{T}: \mathcal{X}\to \mathcal{X}$ satisfying 
\begin{enumerate}
\item $\mathcal{T}\circ\mathcal{T}(\Phi) = \Phi$ 
for any $\Phi\in \mathcal{X}$,  
\item For any $\Phi\in \mathcal{X}$, 
$\mathcal{T}(\Phi)=\Phi$  if $\Phi\in \mathcal{Y}$ and 
$\mathcal{T}(\Phi) \ne \Phi$ if $\Phi\not\in \mathcal{Y}$,      
%$\mathcal{GY}$.   
\item $\mathcal{T}\left(\mathcal{GX}\right)\subset\mathcal{X}$  and 
$\mathcal{T}\left(\mathcal{GX}\right)\not\subset\mathcal{GX}$.    
\end{enumerate} 
Such an involution is called a \textit{fake Legendre involution}.   
For details on fake Legendre involution, see CASE II-1, CASE II-2 in \S 3.     
For a fake Legendre involution $\mathcal{T}: \mathcal{X}\to \mathcal{X}$, 
the property (3) implies $\mathcal{T}\ne \mathcal{L}$.     
 Since the main purpose of this paper is to obtain a characterization of 
the Legendre transform $\mathcal{L}$, we need to introduce 
$\mathcal{GX}$.   
\par  
%We want to obtain a characterization of the Legendre involution on the class 
%$\mathcal{X}=\mathcal{C}\left(\mathcal{F}\right)$.    
%However, in order to avoid 
%fake Legendre involutions (for fake Legendre involutions, 
%see CASE II-1, CASE II-2 in \S 3), 
%Definiotn \ref{definition1} is needed.      
%Thus, Definition \ref{definition1} gives a technical condition 
%for the proof of the main theorem.    
\medskip 
This paper is organized as follows.
In Section 2, preliminaries are given.   
Theorem \ref{theorem2} is proved in Section \ref{section3}.
%%%%%%%%%%%%%%%%%%%%%%%%%%%%%%%%%%%%%%%%%%%%%%% 
%%%%%%%%%%%%%%%%%%%%%%%%%%%%%%%%%%%%%%%%%%%%%%% 
\section{Preliminaries}\label{section2}  
%%%%%%%%%%%%%%%%%%%%%%%%%%%%%%%%%%%%%%%%%%%%%%% 
\subsection{Envelopes created by hyperplane families}
\label{subsection2.1}  
%%%%%%%%%%%%%%%%%%%%%%%%%%%%%%%%%%%%%%%%%%%%%%%   
Given a function $a: \R^n\to \R$ and a mapping 
$\nu: \R^n\to S^n$, for any $x\in \R^n$ the hyperplane 
$H_{\left(\nu(x), a(x)\right)}$ is 
naturally defined as follows.   
\[
H_{\left(\nu(x), a(x)\right)}=\left\{\left.X\in \R^{n+1}\; \right|
\langle X, \nu(x)\rangle=a(x)\right\}.    
\]
\begin{definition}[Definition 1 in \cite{nishimura}]\label{definition3}
{\rm 
A mapping $\varphi: \R^n\to \R^{n+1}$ is called an 
\textit{envelope} created by 
 the hyperplane family 
$\left\{H_{\left(\nu(x), a(x)\right)}\right\}_{x\in \R^n}$ if the following 
two conditions hold.   
\begin{enumerate}
\item $\varphi(x)\in H_{\left(\nu(x), a(x)\right)}$ for any $x\in \R^n$.   
\item $\langle d\varphi_x({\bf v}), \nu(x)\rangle=0$ for any $x\in \R^n$ and 
any ${\bf v}\in T_x\R^n$.    
\end{enumerate}
} 
\end{definition}
\begin{definition}[A strong version of Definition 2 in \cite{nishimura}]
\label{definition4}
{\rm 
The hyperplane family 
$\left\{H_{\left(\nu(x), a(x)\right)}\right\}_{x\in \R^n}$ is said to be 
\textit{creative} if there exist a mapping 
$\Omega: \R^n\to T^*S^n$ with the form 
$\Omega(x)=\left(\nu(x), \omega(x)\right)$ such that 
the equality $da=\omega$ holds for any $x\in \R^n$, 
where $T^*S^n$ is the cotangent bundle of $S^n$.   
}
\end{definition}
\begin{remark}
\begin{enumerate}
\item Notice that in Definition \ref{definition4}, the equality 
$da=\omega$ globally holds, while in Definition 2 of \cite{nishimura} 
the equality $da=\omega$ holds as germs of one-form at $x$ 
for any $x\in \mathbb{R}^n$.   
Thus,Definition \ref{definition4} is apparently 
a strong version of Definition 2 in  
\cite{nishimura}.    
\item Despite (1) above, since the source space of $a, \nu, \Omega$ 
is a convex set $\R^n$, it is possible to obtain the concrete expression 
of the one-form $\omega$ by using a single normal 
cordinate neighborhood 
$\left(V, \Theta_1, \ldots, \Theta_n\right)$ as follows, 
where $V$ is a small neighborhood 
of ${\bf e}_1=\left(1, 0, \ldots, 0\right)\in S^n$.      \\
Notice that if $\nu: \R^n\to S^n$ is a Gauss mapping of a frontal 
$\varphi$, then $-\nu$, too, is a Gauss mapping of $\varphi$.   
Thus, without loss of generality, from the first we may assume 
$\nu(0)\ne -{\bf e}_1=\left(-1, 0, \ldots, 0\right)$.   
Therefore, the Levi-Civita translation 
$\Pi_{\left(\nu(0), {\bf e}_1\right)}: T_{{\bf e}_1}S^n\to T_{\nu(0)}S^n$ 
is uniquely determined.      
Given a non-zero $x_0\in \R^n$, set $x_t=tx_0$ $(t\in [0, 1])$.    
Since the Levi-Civita translation 
$\Pi_{\left(\nu(x), \nu(0)\right)}: T_{\nu(0)}S^n\to T_{\nu(x)}S^n$  
is continuously depending on $x\in \R^n$ and the interval 
$[0,1]$ is connected, 
the Levi-Civita translation 
\[
\Pi_{\left(\nu(x_t), {\bf e}_1\right)}: T_{{\bf e}_1}S^n\to T_{\nu(x_t)}S^n 
\]
is uniquely determined for any $t\in [0,1]$.    
For any $x\in \R^n$ and any $i\in \{1, \ldots, n\}$, set 
\[
b_i(x)  =  \omega(x)\left(\Pi_{\left(\nu(x), {\bf e}_1\right)}
\left(\left(\frac{\partial}{\partial \Theta_i}\right)
_{\mbox{at ${\bf e}_1$}}\right)\right) \;\mbox{ and }\;   
d\theta_i  =  d\left(\Theta_{(x, i)}\circ \nu\right) 
\mbox{ at $x$},   
\]
where $\left(\frac{\partial}{\partial \Theta_{(x, i)}}\right)= 
\Pi_{\left(\nu(x), {\bf e}_1\right)}
\left(\left(\frac{\partial}{\partial \Theta_i}\right)
_{\mbox{at ${\bf e}_1$}}\right)$.    Then, it follows 
\[
\omega(x)  =  \sum_{i=1}^n b_i(x)d\theta_i \; \mbox{ and }\;  
\nu(x)  =  \exp\left(\theta_1(x), 
\ldots, \theta_n(x)\right), 
\]
where $\exp: T_{{\bf e}_1}S^n\to S^n$ is the exponential mapping 
at ${\bf e}_1$.    
\end{enumerate}
\end{remark}
\begin{theorem}[(a) of Thorem 1 in \cite{nishimura}]\label{theorem3}
The hylerplane family 
$\left\{H_{\left(\nu(x), a(x)\right)}\right\}_{x\in \R^n}$ creates an envelope 
if and only if it is creative.   
\end{theorem}
\begin{theorem}[(b) of Theorem 1 in \cite{nishimura}]\label{theorem4}
Suppose that the hylerplane family 
$\left\{H_{\left(\nu(x), a(x)\right)}\right\}_{x\in \R^n}$ creates an envelope 
$\varphi: \mathbb{R}^n\to \mathbb{R}^{n+1}$.   
Then, for any $x\in \R^n$, under the canonical identifications 
$T^*_{\nu(x)}S^n\cong T_{\nu(x)}S^n \subset T_{\nu(x)}\mathbb{R}^{n+1}
\cong \R^{n+1}$, the $(n+1)$-dimensional vector $\varphi(x)$
is represented as follows.
\[
\varphi(x) = \omega(x) + a(x)\nu(x), 
%\exp\left(\theta_1(x), \ldots, \theta_n(x)\right),  
\]
where the $(n+1)$-dimensional vector 
$\omega(x)$ is identified with the corresponding $n$-dimensional 
cotangent vector $\omega(x)=\sum_{i=1}^nb_i(x)d\theta_i$ 
under these identifications.
\end{theorem}
\begin{example}\label{example1}
Let $\varphi: \R\to \R^2$, $\varphi(x)=\left(x^2, x^5\right)$ 
be a frontal given in \S 1.   Let 
$\nu: \R\to S^1$ be the mapping defined by 
$\nu(x)=\left(\frac{-5x^3}{\sqrt{25x^6+4}}, \frac{2}{\sqrt{25x^6+4}}\right)$.    
It is clear that $\nu$ is a Gauss mapping of $\varphi$.    
Set $a(x)=\langle \varphi(x), \nu(x)\rangle=\frac{-3x^5}{\sqrt{25x^6+4}}$.    
%$\left(\cos \theta(x), \sin\theta(x)\right)=\nu(x)$.    
Then, the straight line 
\[
H_{\left(\nu(x), a(x)\right)}=
\left\{\left(X, Y\right)\in \R^2\; \left|\; 
\langle (X, Y), \nu(x)\rangle =a(x)\right.\right\}
\]
is the affine tangent line to the image of $\varphi$ if $x\ne 0$ and 
$H_{\left(a(0), \nu(0)\right)}$ is well-defined.   
From the construction, the given frontal must be an 
envelope of the family 
$\left\{H_{\left(\nu(x), a(x)\right)}\right\}_{x\in \R}$.    
The following is the confirmation 
of this fact by using Theorem \ref{theorem3} and 
Theorem \ref{theorem4}.     
Set $\left(\cos \theta(x), \sin\theta(x)\right)=\nu(x)$.    
Calculations show the following.   
\[
\frac{da}{dx}(x)=\frac{-150x^{10}-60x^4}{\left(25x^6+4\right)^{\frac{3}{2}}}, 
\quad 
\frac{d\theta}{dx}(x)=\frac{30x^2}{25x^6+4}.    
\]
Since $\frac{da}{dx}(x)$ can be divided by $\frac{d\theta}{dx}(x)$, by Theorem 
\ref{theorem3}, the line family 
$\left\{H_{\left(\nu(x), a(x)\right)}\right\}_{x\in \R}$ creates an envelope.    
Set 
\[
b(x)=\frac{\frac{da}{dx}(x)}{\frac{d\theta}{dx}(x)}=
\frac{-x^2\left(5x^6+2\right)}{\sqrt{25x^6+4}}.   
\]
By applying Theorem \ref{theorem4}, it is seen that 
the image of $\varphi$ is actually an envelope of 
$\left\{H_{\left(\nu(x), a(x)\right)}\right\}_{x\in \R}$ as follows.   
\begin{eqnarray*}
{} & { } & 
a(x)\left(\cos\theta(x), \sin\theta(x)\right)+ 
b(x)\left(-\sin\theta(x), \cos\theta(x)\right) \\ 
{ } & = & 
\frac{-3x^5}{\sqrt{25x^6+4}}
\left(\frac{-5x^3}{\sqrt{25x^6+4}}, \frac{2}{\sqrt{25x^6+4}}\right)
+\frac{-2x^2\left(5x^6+2\right)}{\sqrt{25x^6+4}}
\left(\frac{-2}{\sqrt{25x^6+4}}, \frac{-5x^3}{\sqrt{25x^6+4}}\right) \\ 
{ } & = & 
\frac{1}{{25x^6+4}}
\left(
15x^8+2x^2\left(5x^6+2\right), -6x^5+5x^5\left(5x^6+2\right)
\right) \\ 
{ } & = & 
\frac{1}{{25x^6+4}}
\left(
15x^8+10x^8+4x^2, -6x^5+25x^{11}+10x^5\right) \\ 
{ } & = & 
\frac{1}{{25x^6+4}}
\left(
x^2\left(25x^6+4\right), x^5\left(25x^6+4\right)\right) \\ 
{ } & = & 
\left(x^2, x^5\right).     
\end{eqnarray*}   
\end{example}
\begin{remark}\label{remark2}
\begin{enumerate}
\item Notice that in Example \ref{example1}, 
$\nu(x)=\left(\cos\theta(x), \sin\theta(x)\right)=
\left(\frac{-5x^3}{\sqrt{25x^6+4}}, \frac{2}{\sqrt{25x^6+4}}\right)$ and 
$b(x)=\frac{-x^2\left(5x^6+2\right)}{\sqrt{25x^6+4}}$,  
%is    
%$\theta(x)=\frac{5}{2}x^3+\mbox{higher}$ 
%(resp., $b(x)=-x^2+\mbox{higher}$).     
%This simple observation shows that 
From these concrete forms, it is seen that there are no sections  
$\xi: S^1\to T^*S^1$ such that $\Omega=\xi\circ \nu$.    
In other words, there are no functions $\beta: S^1\to \R$ such that 
$b(x)=\beta\circ \nu(x)$.    
As a conclusion of this simple observation, 
it is seen that the function $b$ can not be obtained by using a contact 
structure of $T^*S^1\times \R$.    
\item Despite (1) of this remark, the function $b$ is directly obtained 
by the creative 
condition $da=\omega$ in Definition \ref{definition4}.   
%creative condition $(*)$.   
This is because Theorem \ref{theorem3} and 
Theorem \ref{theorem4} are proved without using 
any information of the whole of cotangent space $T^*S^n$.   
%These are proved by using only 
%information of cotangent vector along the Gauss mapping 
%$\nu: \R^n\to S^n$.   
Theorem \ref{theorem3} and 
Theorem \ref{theorem4} are proved by using the anti-orthotomic method 
developed in \cite{janeczkonishimura} and the anti-orthotomic method  
needs only information of cotangent vectors along the Gauss mapping 
$\nu: \R^n\to S^n$.
\end{enumerate}
\end{remark}  
%\begin{theorem}[Theorem 2 in \cite{nishimura}]\label{theorem5}
%Suppose that the hylerplane family 
%$\left\{H_{\left(\nu(x), a(x)\right)}\right\}_{x\in \R^n}$ creates an envelope 
%$\varphi: \mathbb{R}^n\to \mathbb{R}^{n+1}$.     
%Then, the created envelope is unique if and only if 
%the set consisting of regular points of $\nu$ is dense in $\R^n$.    
%\end{theorem}
%%%%%%%%%%%%%%%%%%%%%%%%%%%%%%%%%%%%%%%%%%%%%%% 
\subsection{Natural correspondence 
$\mathcal{C}: \mathcal{F}\to \mathcal{X}$}\label{subsection2.2}
%%%%%%%%%%%%%%%%%%%%%%%%%%%%%%%%%%%%%%%%%%%%%%% 
Recall that 
\begin{eqnarray*}
\mathcal{F} & = & \left\{frontls\right\} \\  
\mathcal{X}  & =&   
\left\{\Phi=\left(\left(\theta_1, \ldots, \theta_n\right), a, 
\left(b_1, \ldots, b_n\right)\right): \mathbb{R}^n\to \mathbb{R}^{2n+1} 
\; \left|\; 
\begin{array}{c}
\theta_i, a, b_i: \mathbb{R}^n\to \mathbb{R}\;\; C^\infty\;\; 
(\forall i\in \{1, \ldots, n\}), \\ 
(*) \mbox{ in \S 1 is satisfied}. \\ 
%Reg(\theta), Reg(b) \mbox{ are dense in $\mathbb{R}^n$}
\end{array}
\right.\right\}.    
\end{eqnarray*}
By using Subsection 2.1, a bijective mapping   
$\mathcal{C}: \mathcal{F}\to \mathcal{X}$ 
%and 
%$\mathcal{C}_2: \mathcal{GX}\to \mathcal{GF}$ 
is naturally constructed as follows. \\ 
\par 
\noindent 
 %%%%%%%%%%%%%%%%%%%%%%%%%%%%%%%%%%%%%%%%%%%%
\underline{\bf Construction of 
$\mathcal{C}: \mathcal{F}\to \mathcal{X}$.
} \\ 
\par   
Let $\varphi$ be an element of $\mathcal{F}$.    
Let $\nu: \R^n\to S^n$ be a Gauss mapping of 
$\varphi$.     Without loss of generality, 
we may assume $\nu(0)\ne -{\bf e}_1=(-1, 0, \ldots, 0)\in S^n$.   
Set $a(x)=\langle \varphi(x), \nu(x)\rangle$.    
Since the given frontal $\varphi$ is an envelope of 
the hyperplane family 
$\left\{H_{\left(\nu(x), a(x)\right)}\right\}_{x\in \mathbb{R}^n}$, 
by Theorem \ref{theorem3}, 
$\left\{H_{\left(\nu(x), a(x)\right)}\right\}_{x\in \mathbb{R}^n}$ is 
creative.   Thus, there exists a mapping 
$\Omega: \mathbb{R}^n\to T^*S^n$ with the form 
$\Omega(x)=\left(\nu(x), \omega(x)\right)$ such that 
$da=\omega$.     
For any $x\in \R^n$, set $x_t=tx\in \mathbb{R}^n$ 
$(t\in [0,1])$.   
Then, since $[0,1]$ is connected, by using the 
equality $\nu(x_t)=
\exp\left(\theta_1(x_t), \ldots, \theta_n(x_t)\right)$ 
where $\exp: T_{{\bf e}_1}S^n\to S^n$ is the exponential mapping, 
for any $t\in [0,1]$ 
the vector $\left(\theta_1(x_t), \ldots, \theta_n(x_t)\right)\in 
T_{{\bf e}_1}S^n$ 
is uniquely determined.    
%Moreover, since $\varphi$ is generic, 
%$Reg(\theta)$ is dense in $\R^n$.     
Thus, by 
Theorem \ref{theorem4} and Theorem \ref{theorem5}, 
the mapping $\left(b_1, \ldots, b_n\right): \R^n\to \R^n$, too, is uniquely 
determined.   Under the canonical 
identification $T_{{\bf e}_1}S^n\cong \R^n$, 
we set 
\[
\Phi=\left(\left(\theta_1, \ldots, \theta_n\right), a, 
\left(b_1, \ldots, b_n\right)\right): \R^n\to \R^{2n+1}.    
\]
%Since $\varphi$ is generic, it follows that $Reg(b)$ is dense in $\mathbb{R}^n$ 
%as well.   Moreover, 
Then, for the $\Phi$ the creative condition $(*)$ is satisfied 
since the equality $da=\omega$ holds  Hence, the constructed 
$\Phi$ belongs to $\mathcal{X}$.     
Thus, the definition   
\[
\mathcal{C}(\varphi)=\Phi.   
\] 
is well-defined.   
\begin{lemma}\label{lemma1}
The constructed mapping $\mathcal{C}: \mathcal{F}\to \mathcal{X}$ is 
injective.   
\end{lemma}
\proof     
Let $\varphi_1, \varphi_2$ be two elements of $\mathcal{F}$.   
Let $\nu_1, \nu_2: \R^n\to S^n$ be Gauss mappings of $\varphi_1, 
\varphi_2$ respectively.    
Without loss of generality, we may assume that 
$\nu_j(0)\ne -{\bf e}_1=(-1, 0, \ldots, 0)$ for any $j\in \{1, 2\}$.     
Set $a_j(x)=\langle \varphi_j(x), \nu_j(x)\rangle$ for any $j\in \{1, 2\}$.    
Suppose that $\mathcal{C}(\varphi_1)=\mathcal{C}(\varphi_2)$.    
Then, from the above construction, it follows that $a_1=a_2$ and 
$\nu_1=\nu_2$. For any $j\in \{1, 2\}$, 
let $\Omega_j: \R^n\to T^*S^n$ be the mapping having the form 
$\Omega_j(x)=\left(\nu_j(x), \omega_j(x)\right)$ satisfying 
$da_j=\omega_j$.    Then, since $a_1=a_2$, %from the above construction, 
it follows $\omega_1=\omega_2$.    Therefore, by Theorem \ref{theorem4}, 
$\varphi_1(x)=
\omega_1(x)+a_1(x)\nu_1(x)=\omega_2(x)+a_2(x)\nu(x)=\varphi_2(x)$.
\hfill $\Box$
\begin{lemma}\label{lemma2}
The constructed mapping $\mathcal{C}: \mathcal{F}\to \mathcal{X}$ is 
surjective.   
\end{lemma}
\proof    
Let $\Phi=\left(\left(\theta_1, \ldots, \theta_n\right), a, 
\left(b_1, \ldots, b_n\right)\right): \R^n\to \R^{2n+1}$ be an element of 
$\mathcal{X}$.     %Set ${\bf e}_1=(1, 0, \ldots, 0)\in S^n$.    
Consider the $(n+1)$-dimensional unit vector 
${\bf v}_0=\exp\left(\theta_1(0), \ldots, \theta_n(0)\right)$ where 
$\exp: T_{{\bf e}_1}S^n\to S^n$ is the exponential mapping.      
If ${\bf v}_0\ne -{\bf e}_1$, then 
define the mapping $\nu: \R^n\to S^n$ by 
$\nu(x)=
\exp\left(\theta_1(x), \ldots, \theta_n(x)\right)$.   
If ${\bf v}_0= -{\bf e}_1$, then 
define the mapping $\nu: \R^n\to S^n$ by 
$\nu(x)=
\exp\left(\theta_1(x)-\theta(0), \ldots, \theta_n(x)-\theta(0)\right)$. 
%Then, it follows $\nu(0)={\bf e}_1$.   
Define the $n$-dimensional cotangent vector 
$\omega(x)\in T_{\nu(x)}^*S^n$ 
and the mapping $\Omega: \R^n\to T^*S^n$ 
by $\omega(x)=\sum_{i=1}^n b_i(x)d\theta_i$ and 
$\Omega(x)=\left(\nu(x), \omega(x)\right)$ respectively.   
Then, since 
$\Phi\in \mathcal{X}$, it follows $da=\omega$.   
% and 
%that $Reg(\theta)$ is dense.    
Thus, by Theorem \ref{theorem3}, % and Theorem \ref{theorem5}, 
the hyperplane family 
$\left\{H_{(\nu(x), a(x))}\right\}_{x\in \R^n}$ creates an envelope.   
Set $\varphi(x)=\omega(x)+a(x)\nu(x)$.    Then, 
by Theorem \ref{theorem4}, the created envelope of 
$\left\{H_{(\nu(x), a(x))}\right\}_{x\in \R^n}$ must be $\varphi$.    
Since $\Phi$ is an element of $\mathcal{X}$ and 
any envelope of hyperplane family is a frontal, 
$\varphi$ is an element of $\mathcal{F}$.   
From the construction, it follows $\mathcal{C}(\varphi)=\Phi$.    
\hfill $\Box$    \\
 
%%%%%%%%%%%%%%%%%%%%%%%%%%%%%%%%%%%%%%%%%%%%
%\underline{\bf Construction of $\mathcal{C}_2: 
%\mathcal{GX}\to \mathcal{GF}$.
%} \\ 
%\par    
%Let $\Phi$ be an element of $\mathcal{GX}$.    
%Then, $\Phi$ has the form 
%\[
%\Phi=\left(\left(\theta_1, \ldots, \theta_n\right), a, 
%\left(b_1, \ldots, b_n\right)\right): \R^n\to \R^{2n+1}
%\]  
%satisfying the following three conditions.   
%\begin{enumerate}
%\item[(a)] $da=\sum_{i=1}^n b_id\theta_i$, 
%\item[(b)] $Reg(\theta)$ is dense in $\R^n$, 
%\item[(c)] $Reg(b)$ is dense in $\R^n$.   
%\end{enumerate}
%Under the identification 
%as $\omega=\sum_{i=1}^nb_id\theta_i$ 
%%%%%%%%%%%%%%%%%%%%%%%%%%%%%%%%%%%%%%%%%%%%  
\section{Proof of Theorem \ref{theorem2}}\label{section3}
%%%%%%%%%%%%%%%%%%%%%%%%%%%%%%%%%%%%%%%%%%%%%%% 
Set 
\begin{eqnarray*}
F  & = & \left(F_1, F_2, F_3\right) : \mathbb{R}^{2n+1}\to \mathbb{R}^{2n+1},  \\ 
F_1 & = & \left(F_{11}, \ldots, F_{1n}\right): 
\mathbb{R}^{2n+1}\to \mathbb{R}^n,  \\ 
F_2 & : & \mathbb{R}^{2n+1}\to \mathbb{R},  \\
F_3 & = & \left(F_{31}, \ldots, F_{3n}\right): \mathbb{R}^{2n+1}\to \mathbb{R}^n 
\\ 
\end{eqnarray*}
and  
\begin{eqnarray*}
\widetilde{\theta}_i & = & F_{1i}\circ
\left(\theta_1, \ldots, \theta_n, a, b_1\ldots, b_n\right)
\quad (\mbox{for any }i\in \{1, \ldots, n\}),  \\
\widetilde{a} & = & 
F_{2}\circ\left(\theta_1, \ldots, \theta_n, a, b_1\ldots, b_n\right) \\ 
\widetilde{b}_i & = & F_{3i}\circ
\left(\theta_1, \ldots, \theta_n, a, b_1\ldots, b_n\right)
\quad (\mbox{for any }i\in \{1, \ldots, n\}).   
\end{eqnarray*} 
For any $i\in \{1, \ldots, n\}$, the coefficient of the term 
$X_1^{j_1}\cdots X_n^{j_n}Y^kZ_1^{\ell_1}\cdots Z_n^{\ell_n}$ of the analytic 
function $F_{1i}$ (resp., $F_{3i}$) is denoted by 
$p_{i, (j_1, \ldots, j_n,k,\ell_1, \ldots, \ell_n)}$ 
(resp., $q_{i, (j_1, \ldots, j_n,k,\ell_1, \ldots, \ell_n)}$).   
Namely, 
\begin{eqnarray*}
F_{1i}\left(X_1, \ldots, X_n, Y, Z_1, \ldots, Z_n\right) 
& = & 
\sum_{0\le j_m, k, \ell_m<\infty} 
p_{i, (j_1, \ldots, j_n,k,\ell_1, \ldots, \ell_n)}
X_1^{j_1}\cdots X_n^{j_n}Y^kZ_1^{\ell_1}\cdots Z_n^{\ell_n} \\ 
F_{3i}\left(X_1, \ldots, X_n, Y, Z_1, \ldots, Z_n\right) 
& = & 
\sum_{0\le j_m, k, \ell_m<\infty} 
q_{i, (j_1, \ldots, j_n,k,\ell_1, \ldots, \ell_n)}
X_1^{j_1}\cdots X_n^{j_n}Y^kZ_1^{\ell_1}\cdots Z_n^{\ell_n}
\end{eqnarray*} 
%we set 
%\[
%F_i(X, Y, Z)=
%\sum_{0\le j, k, \ell <\infty} p_{i, (j, k, \ell)}X^jY^kZ^\ell, 
%\]
%where $p_{i, (j, k, \ell)}\in \mathbb{R}$ and 
We define 
$X_i^0=Y^0=Z_i^0=1$ for any $i\in \{1, \ldots, n\}$.    
%%%%%%%%%%%%%%%%%%%%%%%%%%%%%%%%%%%%%%%%%%%%% 
For any $i\in \{1, \ldots, n\}$, let ${\bf e}_i$ be the $i$-th row vector of the 
$n$ by $n$ unit matrix and ${\bf 0}$ is the zero vector in the vector space 
$\mathbb{R}^n$.    
Then, it is possible to abbreviate 
coefficients $p_{i, (1, 0, \ldots, 0)}$, $q_{i, (0, \ldots, 0, 1)}$ etc. as 
$p_{i, ({\bf e}_1, 0, {\bf 0})}$,  $q_{i, ({\bf 0}, 0, {\bf e}_n)}$ etc. 
respectively, and 
we adopt such abbreviations.   
\par 
The fixed point set assumption 
(the assumption (2) of Theorem \ref{theorem2}) implies the following:  
\begin{lemma}\label{lemma3.1}
The following two hold for any $i\in \{1, \ldots, n\}$.    
\begin{eqnarray*}
F_{1i}(X_1, \ldots, X_n,Y, Z_1, \ldots, Z_n) & = & 
p_{i, ({\bf e}_i,0,{\bf 0})}X_i+p_{i, ({\bf 0},0,{\bf e}_i)}Z_i, \\  
%\quad \left(p_{1, (1,0,0)}+p_{1, (0,0,1)}=1\right), \\ 
F_{3i}(X_1, \ldots, X_n, Y, Z_1, \ldots, Z_n) & = & 
q_{i, ({\bf e}_i,0,{\bf 0})}X_i+q_{i, ({\bf 0},0,{\bf e}_i)}Z_i.    \\ 
%F_3(X, Y, Z) & = & p_{3, (1,0,0)}X+p_{3, (0,0,1)}Z 
%\quad \left(p_{3, (1,0,0)}+p_{3, (0,0,1)}=1\right).    
\end{eqnarray*}
\end{lemma}
The assumptions (1), (2) of Theorem \ref{theorem2} imply the following:   
%involution assumption 
%(the assumption (1) in Theorem \ref{theorem2}) implies the following: 
\begin{lemma}\label{lemma3.2}%\quad %\\ 
The following six hold for any $i\in \{1, \ldots, n\}$.    
\begin{enumerate}
\item $p_{i, ({\bf e}_i, 0, {\bf 0})}+p_{i, ({\bf 0}, 0, {\bf e}_i)}=1$,  \\ 
\item $q_{i, ({\bf e}_i, 0, {\bf 0})}+q_{i, ({\bf 0}, 0, {\bf e}_i)}=1$,  \\ 
\item $p_{i, ({\bf e}_i, 0, {\bf 0})}^2
+p_{i, ({\bf 0},0,{\bf e}_i)}q_{i, ({\bf e}_i,0,{\bf 0})}=1$, \\ 
\item $p_{i, ({\bf e}_i, 0, {\bf 0})}p_{i, ({\bf 0}, 0, {\bf e}_i)}
+p_{i, ({\bf 0}, 0, {\bf e}_i)}q_{i, ({\bf 0}, 0, {\bf e}i)}=0$, \\ 
\item $q_{i, ({\bf e}_i, 0, {\bf 0})}p_{i, ({\bf e}_i, 0, {\bf 0})}
+q_{i, ({\bf 0}, 0, {\bf e}_i)}q_{i, ({\bf e}_i, 0, {\bf 0})}=0$, \\
\item $q_{i, ({\bf e}_i, 0, {\bf 0})}p_{i, ({\bf 0}, 0, {\bf e}_i)}
+q_{i, ({\bf 0},0, {\bf e}_i)}^2=1$. 
\end{enumerate} 
\end{lemma}
It is easy to confirm that in Lemma \ref{lemma3.2}, 
under the equations (1) and (2), the equation (3) 
(resp., (5)) is 
equivalent to the equation (4) (resp., (6)).   
%under the equation 
%(1), and the equation (5) is 
%equivalent to the equation (6) under the equation 
%(2). 
%$p_{3, (1,0,0)}+p_{3, (0,0,1)}=1$.      
From now on, we concentrate on considering  
 (1), (2), (3) and (5) for any $i\in \{1, \ldots, n\}$.   
Namely, for any $i\in \{1, \ldots, n\}$, we %first 
seek the solution set 
of 
% $p_{1, (1,0,0)}$, $p_{1, (0,0,1)}$, $p_{3, (1,0,0)}$ and $p_{3, (0,0,1)}$, 
the following system of quadratic equations:   
\[
\left\{
\begin{array}{ccc}
p_{i, ({\bf e}_i, 0, {\bf 0})}+p_{i, ({\bf 0}, 0, {\bf e}_i)} & = & 1, \\ 
q_{i, ({\bf e}_i, 0, {\bf 0})}+q_{i, ({\bf 0}, 0, {\bf e}_i)} & = & 1,  \\ 
p_{i, ({\bf e}_i, 0, {\bf 0})}^2
+p_{i, ({\bf 0},0,{\bf e}_i)}q_{i, ({\bf e}_i,0,{\bf 0})} & = & 1,  \\ 
%p_{1, (1,0,0)}^2+p_{1, (0,0,1)}p_{3, (1,0,0)} & = & 1,  \\ 
q_{i, ({\bf e}_i, 0, {\bf 0})}p_{i, ({\bf e}_i, 0, {\bf 0})}
+q_{i, ({\bf 0}, 0, {\bf e}_i)}q_{i, ({\bf e}_i, 0, {\bf 0})} & = & 0. 
%p_{3, (1,0,0)}p_{1, (1,0,0)}+p_{3, (0,0,1)}p_{3, (1,0,0)} & = & 0.   
\end{array}
\right.
\]
It is not difficult to solve this system of quadratic equations completely.   
The solution set of this system of equations is the disjoint union 
of the following two sets.   
\begin{center}
%\begin{itemize}
%\item[{\bf (a)}]\quad 
%\begin{eqnarray*}
%\left\{\left(p_{i, ({\bf e}_i, 0, {\bf 0})}, 
%p_{i, ({\bf 0}, 0, {\bf e}_i)}, 
%q_{i, ({\bf e}_i, 0, {\bf 0})}, 
%q_{i, ({\bf 0}, 0, {\bf e}_i)}\right)\in\mathbb{R}^4\right\}
%& = & 
\[ 
\left\{\left(1,0,0,1\right)\right\}, \quad  
%\item[{\bf (b)}]\quad 
%\left\{\left(p_{i, ({\bf e}_i, 0, {\bf 0})}, 
%p_{i, ({\bf 0}, 0, {\bf e}_i)}, 
%q_{i, ({\bf e}_i, 0, {\bf 0})}, 
%q_{i, ({\bf 0}, 0, {\bf e}_i)}\right)\in\mathbb{R}^4\right\} 
%& = & 
\left\{\left.\left(t_i-1, 2-t_i, t_i, 1-t_i\right)\in \mathbb{R}^4\; \right|\; 
t_i\in \mathbb{R}\right\}.
\]   
%\end{eqnarray*} 
%\end{itemize}
\end{center} 
%In the case {\bf (i)}, it follows $\widetilde{\theta}_i=\theta_i$ and 
%$\widetilde{b}_i=b_i$.     On the other hand, in the case of {\bf (ii)}, 
%we have that $\widetilde{\theta}_i=\left(t_i-1\right)\theta_i+
%\left(2-t_i\right)b_i$ and $\widetilde{b}_i=t_i\theta_i+
%\left(1-t_i\right)b_i$.    
\par 
\medskip 
\noindent 
\underline{\bf CASE I. 
$\exists i\in \{1, \ldots, n\}$ such that 
$\left\{\left(p_{i, ({\bf e}_i, 0, {\bf 0})}, 
p_{i, ({\bf 0}, 0, {\bf e}_i)}, 
q_{i, ({\bf e}_i, 0, {\bf 0})}, 
q_{i, ({\bf 0}, 0, {\bf e}_i)}\right)\in\mathbb{R}^4\right\}
=\left\{\left(1,0,0,1\right)\right\}$.} \\ 
\par 
Let $i_0\in \{1, \ldots, n\}$ be the ineger satisfying 
\[
\left\{\left(p_{i_0, ({\bf e}_{i_0}, 0, {\bf 0})}, 
p_{{i_0}, ({\bf 0}, 0, {\bf e}_{i_0})}, 
q_{i_0, ({\bf e}_{i_0}, 0, {\bf 0})}, 
q_{i_0, ({\bf 0}, 0, {\bf e}_{i_0})}\right)\right\}
=\left\{\left(1,0,0,1\right)\right\}.
\]   
For the $i_0$, it follows 
$\widetilde{\theta}_{i_0}=\theta_{i_0}$ and 
$\widetilde{b}_{i_0}=b_{i_0}$.   % for any $i\in \{1, \ldots, n\}$.      
%\par    
Let 
$\Phi=\left(\left(\theta_1, \ldots, \theta_n\right), 
\frac{1}{2}\sum_{i=1}^n\theta_i^2, \left(\theta_1, \ldots, \theta_n\right)\right)$ 
be an element of $\mathcal{GY}$.    For the $\Phi$, set 
\[
\Phi_0=\left({\bf 0}, \theta_0^2, 2\theta_{i_0}{\bf e}_{i_0}\right).   
\]   
Then, it is clear that $\varphi+\lambda\varphi_0\not\in \mathcal{GY}$ 
if $\lambda\ne 0$.      
On the other hand, since $\widetilde{\theta}_{i_0}=\theta_{i_0}$ and 
$\widetilde{b}_{i_0}=b_{i_0}$, it follows 
$\mathcal{T}(\Phi+\lambda\Phi_{0})=\Phi+\lambda\Phi_{0}$ 
for any $\lambda\in \mathbb{R}$.    
Moreover, it is easily seen that the Jacobian determinant of 
the mapping 
$\left(\theta_1, \ldots, \theta_n\right)+2\lambda\theta_{i_0}{\bf e}_{i_0}$ 
is exactly $\left(1+2\lambda\right)$ times the Jacobian determinant 
of $\left(\theta_1, \ldots, \theta_n\right)$ for any $\lambda\in \mathbb{R}$.   
Hence, $\varphi+\lambda\varphi_{0}\in \mathcal{GX}$.    
Therefore, the given transform $\mathcal{T}$ does not satisfy 
the second half of the assumption (2) of Theorem \ref{theorem2}, and 
in this case there are no transforms satisfying the assumptions (1) and (2) of 
Theorem \ref{theorem2}.      
%%%%%%%%%%%%%%%%%%%%%%%%%%%%%%%%%%%%%%%%%%%%%%%    
%It follows 
%\[
%d\widetilde{a}=\sum_{i=1}^n \widetilde{b}_id\widetilde{\theta}_i 
%= \sum_{i=1}^n {b}_id{\theta}_i=da.    
%\]
%Thus, we have 
%\[
%\widetilde{a}=a+c_0 
%\]
%where $c_0$ is a constant real number.   
%By the first half of the assumption (2) of Theorem \ref{theorem2}, 
%it follows $c_0=0$.    Thus, the given transform 
%$\mathcal{T}$ must be the identity 
%transform, which contradicts the second half of the assumption (2) of 
%Theorem \ref{theorem2}.    
\par 
\medskip 
\noindent 
%%%%%%%%%%%%%%%%%%%%%%%%%%%%%%%%%%%%%%%%%%%%
\underline{\bf CASE II. 
$\left\{\left(p_{i, ({\bf e}_i, 0, {\bf 0})}, 
p_{i, ({\bf 0}, 0, {\bf e}_i)}, 
q_{i, ({\bf e}_i, 0, {\bf 0})}, 
q_{i, ({\bf 0}, 0, {\bf e}_i)}\right)\right\}=
\left\{\left.\left(t_i-1, 2-t_i, t_i, 1-t_i\right)\; \right|\; 
t_i\in \mathbb{R}\right\}$ ($\forall i\in \{1, \ldots, n\}$).    
}   \\ 
\par    
In this case, for any $i\in \{1, \ldots, n\}$, it follows 
\[
\widetilde{\theta}_i=\left(t_i-1\right)\theta_i+\left(2-t_i\right)b_i, 
\quad 
\widetilde{b}_i=t_i\theta_i+\left(1-t_i\right)b_i.     
\]
By the creative condition $(*)$, we have  
\begin{eqnarray*}
{ } & { } & d\widetilde{a} \\ 
{ } & = & \sum_{i=1}^n\widetilde{b}_id\widetilde{\theta}_i \\ 
{ } & = & \sum_{i=1}^n\left(t_i\theta_i+\left(1-t_i\right)b_i\right)
d\left(\left(t_i-1\right)\theta_i+\left(2-t_i\right)b_i\right) \\ 
{ } & = & \sum_{i=1}^n\left(t_i\left(t_i-1\right)\theta_id\theta_i
+t_i\left(2-t_i\right)\theta_idb_i
+\left(1-t_i\right)\left(t_i-1\right)b_id\theta_i
+\left(1-t_i\right)\left(2-t_i\right)b_idb_i\right) \\ 
{ } & = & \sum_{i=1}^n\left(
d\left(\frac{t_i\left(t_i-1\right)}{2}\theta_i^2\right)
+t_i\left(2-t_i\right)\left(d\left(\theta_ib_i\right)-b_id\theta_i\right) 
+\left(1-t_i\right)\left(t_i-1\right)b_id\theta_i 
+d\left(\frac{\left(1-t_i\right)\left(2-t_i\right)}{2}b_i^2\right) 
\right) \\ 
{ } & = & 
\sum_{i=1}^n\left(
d\left(\frac{t_i\left(t_i-1\right)}{2}\theta_i^2\right)
+d\left(t_i\left(2-t_i\right)\theta_ib_i\right)-b_id\theta_i
+d\left(\frac{\left(1-t_i\right)\left(2-t_i\right)}{2}b_i^2\right) 
\right) \\ 
{ } & = & 
d\left(\sum_{i=1}^n
\left(
\frac{t_i\left(t_i-1\right)}{2}\theta_i^2
+ t_i\left(2-t_i\right)\theta_ib_i  
+\frac{\left(1-t_i\right)\left(2-t_i\right)}{2}b_i^2
\right)-a
\right).   
\end{eqnarray*}
\par 
\noindent
%\medskip 
Thus, by using the fixed point set assumption (the assumption (2) of Theorem 
\ref{theorem2}) as well, it follows
\[
\widetilde{a}=
\sum_{i=1}^n
\left(
\frac{t_i\left(t_i-1\right)}{2}\theta_i^2
+ t_i\left(2-t_i\right)\theta_ib_i  
+\frac{\left(1-t_i\right)\left(2-t_i\right)}{2}b_i^2
\right)-a.    
\]   
%%%%%%%%%%%%%%%%%%%%%%%%%%%%%%%%%%%%%%%%%% 
Set 
\[
\widetilde{\widetilde{a}} 
=
F_2\circ \left(\widetilde{\theta}_1, \ldots, \widetilde{\theta}_n, 
\widetilde{a}, \widetilde{b}_1, \ldots, \widetilde{b}_n\right).    
\]
%%%%%%%%%%%%%%%%%%%%%%%%%%%%%%%%%% 
Then, just by substituting, we have 
{\small 
\[
\widetilde{\widetilde{a}} 
=  
\sum_{i=1}^n
\left(\widetilde{A}_i+\widetilde{B}_i+\widetilde{C}_i\right) 
-\left(\sum_{i=1}^n\left(A_i+B_i+C_i\right)-a\right), 
\]
where 
\begin{eqnarray*}
\widetilde{A}_i & = & 
\frac{t_i\left(t_i-1\right)}{2}\left(\left(t_i-1\right)\theta_i+
\left(2-t_i\right)b_i\right)^2, \\ 
\widetilde{B}_i & = & 
t_i\left(2-t_i\right)\left(\left(t_i-1\right)\theta_i+\left(2-t_i\right)b_i\right)
\left(t_i\theta_i+\left(1-t_i\right)b_i\right), \\ 
\widetilde{C}_i & = & 
\frac{\left(1-t_i\right)\left(2-t_i\right)}{2}
\left(t_i\theta_i+\left(1-t_i\right)b_i\right)^2,  \\
A_i & = & 
\frac{t_i\left(t_i-1\right)}{2}\theta_i^2 \\ 
B_i & = & 
t_i\left(2-t_i\right)\theta_ib_i \\ 
C_i & = & 
\frac{\left(1-t_i\right)\left(2-t_i\right)}{2}b_i^2.   
\end{eqnarray*}
}
The function $\widetilde{\widetilde{a}}$ may be represented 
as follows:  
\[
\widetilde{\widetilde{a}} = 
\sum_{i=1}^n
\left(
P_{\theta_i^2}(t_i)\theta_i^2
+P_{\theta_ib_i}(t_i)\theta_ib_i 
+P_{b_i^2}(t_i)b_i^2
\right) +a, 
\] 
where $P_{\theta_i^2}(t_i), P_{\theta_ib_i}(t_i), P_{b_i^2}(t_i)$ are  
apparently quartic polynomials with respect to $t_i$.    
Elementary calculations yield the simplest form of them as follows. 
\begin{eqnarray*}
P_{\theta_i^2}(t_i) & = & \frac{t_i\left(t_i-1\right)^3}{2} 
+ t_i^2\left(2-t_i\right)\left(t_i-1\right) 
-\frac{t_i\left(t_i-1\right)}{2}
+\frac{\left(1-t_i\right)\left(2-t_i\right)t_i^2}{2} \equiv 0, \\ 
P_{\theta_ib_i}(t_i) & = & t_i\left(t_i-1\right)^2\left(2-t_i\right)  
- 2t_i\left(2-t_i\right)\left(t_i-1\right)^2  
+\left(1-t_i\right)^2t_i\left(2-t_i\right)\equiv 0,    \\  
P_{b_i^2}(t_i) & = & \frac{t_i\left(t_i-1\right)\left(2-t_i\right)^2}{2} 
+ t_i\left(2-t_i\right)^2\left(1-t_i\right) 
-\frac{\left(1-t_i\right)\left(2-t_i\right)}{2}
+\frac{\left(1-t_i\right)^3\left(2-t_i\right)}{2} \\ 
{ } & \equiv & \frac{1}{2}t_i\left(t_i-1\right)\left(2-t_i\right)^2.  \\ 
\end{eqnarray*}   
%\vspace{2cm}
On the oher hand, by the involutive assumption 
(the assumption (1) of Theorem \ref{theorem2}), it follows 
\[
\widetilde{\widetilde{a}} 
= a,    
\] 
which is equivalent to assert that for each 
$\i\in \{1, \ldots, n\}$, 
$t_i$ must be a solution of the equation $P_{b_i^2}(t_i)=0$.    
Therefore, for each $i\in \{1, \ldots, n\}$, 
$t_i$ must be $0$ or $1$ or $2$.   \\ 
%contained 
%the following set.    %system of quartic equations: 
%\[
%\left\{
%\begin{array}{ccc}
%P_{\theta_i^2}(t_i) & = & 0 \\ 
%P_{\theta_ib_i}(t_i) & = & 0 \\
%P_{b_i^2}(t_i) & = & 0.   
%\end{array}
%\right.
%\]
%It is easily seen that for any $i\in \{1, \ldots, n\}$ 
%\[
%\left\{t_i\in \mathbb{R}\; \left|\; P_{\theta_i^2}(t_i)=P_{\theta_ib_i}(t_i)=
%P_{b_i^2}(t_i)=0\right.\right\}= \{1, 2\}.   
%\]
%%%%%%%%%%%%%%%%%%%%%%%%%%%%%%%%%%%%% 
\par 
\smallskip 
\noindent 
%%%%%%%%%%%%%%%%%%%%%%%%%%%%%%%%%%%%%%%%%%%%
\underline{\bf CASE II-1. 
$\exists i\in \{1, \ldots, n\}$ such that $t_i=2$.
} \\ 
\par   
%%%%%%%%%%%%%%%%%%%%%%%%%%%%%%%%%%%%%%%%%%%% 
Let $i_0\in \{1, \ldots, n\}$ be the ineger satisfying 
$t_{i_0}=2$.  
For the $i_0$, it follows 
$\widetilde{\theta}_{i_0}=\theta_{i_0}$ and 
$\widetilde{b}_{i_0}=2\theta_{i_0}-b_{i_0}$.   % for any $i\in \{1, \ldots, n\}$.    
For a given $C^\infty$ mapping 
$\Phi_1=\left(\theta_1, \ldots, \theta_n\right): \mathbb{R}^n\to \mathbb{R}^n$ 
such that the set of regular points of it is dense, 
set 
\[
\Phi=\left(\Phi_1, \sum_{i=1}^n \theta_i^2, 2\Phi_1\right): 
\mathbb{R}^n\to \mathbb{R}^{2n+1}.    
\]
Then, $\Phi$ satisfies the creative condition $(*)$ and thus 
$\Phi\in \mathcal{GX}$.  Moreover, since 
$\widetilde{\theta}_{i_0}=\theta_{i_0}$ and 
$\widetilde{b}_{i_0}=2\theta_{i_0}-b_{i_0}$, it follows 
\[
\mathcal{T}(\Phi) = 
\left(\left(\widetilde{\theta}_1, \ldots, \widetilde{\theta}_n\right), 
\widetilde{a}, 
\left(\widetilde{b}_1, \ldots, \widetilde{b}_n\right)
-\widetilde{b}_{i_0}{\bf e}_{i_0}\right).  
\]  
Notice that there are no regular points of the mapping 
$\left(\widetilde{b}_1, \ldots, \widetilde{b}_n\right)
-\widetilde{b}_{i_0}{\bf e}_{i_0}$.      Thus $\mathcal{T}$ 
is 
not an involution  
of $\mathcal{GX}$, it is a fake Legendre involution  
defined in \S 1.   Therefore, in this case, 
there are no transforms $\mathcal{GX}\to \mathcal{GX}$ satisfying 
the assumptions (1) and (2) of Theorem \ref{theorem2}.   \\   
%In this case, $\mathcal{T}$ is called a \textit{fake Legendre involution}.   \\ 
\par 
\smallskip 
\noindent 
%%%%%%%%%%%%%%%%%%%%%%%%%%%%%%%%%%%%%%%%%%%%
\underline{\bf CASE II-2. 
$\exists i\in \{1, \ldots, n\}$ such that $t_i=0$.
} \\ 
\par   
Let $i_0\in \{1, \ldots, n\}$ be the ineger satisfying 
$t_{i_0}=0$.  
For the $i_0$, it follows 
$\widetilde{\theta}_{i_0}=-\theta_{i_0}+2b_i$ and 
$\widetilde{b}_{i_0}=b_{i_0}$.  
From the forms of $\widetilde{\theta}_{i_0}$ and 
$\widetilde{b}_{i_0}$, it is clear that $\mathcal{T}$ 
is a fake Legendre involution.   
Therefore, in this case again, 
there are no transforms $\mathcal{GX}\to \mathcal{GX}$ satisfying 
the assumptions (1) and (2) of Theorem \ref{theorem2}.   \\   
%the strategy of the proof given in 
%the Case II-1 works well in this case as well.   
%In this case as well, $\mathcal{T}$ is called a \textit{fake Legendre involution}.  
% \\ 
%if we take $\Phi$ 
\par 
\smallskip 
\noindent 
%%%%%%%%%%%%%%%%%%%%%%%%%%%%%%%%%%%%%%%%%%%% 
%%%%%%%%%%%%%%%%%%%%%%%%%%%%%%%%%%%%%%%%%%%%
\underline{\bf CASE II-3. 
$t_i=1$ for $\forall i\in \{1, \ldots, n\}$.     
} \\ 
\par   
%%%%%%%%%%%%%%%%%%%%%%%%%%%%%%%%%%%%%%%%%%%% 
In this case, it follows 
$\widetilde{\theta}_{i}=b_{i}$ and 
$\widetilde{b}_{i}=\theta_{i}$ for any $i\in \{1, \ldots, n\}$.   
By the creative condition $(*)$, we have 
\[
d\widetilde{a}=\sum_{i=1}^n\widetilde{b}_i d\widetilde{\theta}_i 
=\sum_{i=1}^n\theta_idb_i 
=d\left(\sum_{i=1}^nb_i\theta_i\right)-\sum_{i=1}^nb_id\theta_i 
=d\left(\sum_{i=1}^nb_i\theta_i\right)-da  
=d\left(\sum_{i=1}^nb_i\theta_i-a\right).    
\]
Therefore, the function $\widetilde{a}$ must have the form 
\[ 
\widetilde{a}=\sum_{i=1}^nb_i\theta_i-a+c_0
\]
where $c_0$ is a constant.   
By the first half of the assumption (2) of Theorem \ref{theorem2}, 
it follows $c_0=0$.     
Therefore, in this case, $\mathcal{T}$ must be the Legendre involution 
$\mathcal{L}$.    This completes the proof.    
\hfill $\Box$
%\vspace{5cm}
%%%%%%%%%%%%%%%%%%%%%%%%%%%%%%%%%%%%%%%%%%%%%%%%%% 
\par 
\bigskip 
As a byproduct of the proof of Theorem \ref{theorem2} given in this section, 
we have a natural complexification of Theorem \ref{theorem2}.    
%Namely, 
%by the proof given in this section, the following corollary is simultaneously 
%obtained.    
%\par 
Set 
\[ 
\mathcal{GX}_{\C}  =  
\left\{\Phi=\left(\left(\theta_1, \ldots, \theta_n\right), a, 
\left(b_1, \ldots, b_n\right)\right): \mathbb{C}^n\to \mathbb{C}^{2n+1} 
\; \mbox{holomorphic} \left|\!%\; 
\begin{array}{c}
%\theta_i, a, b_i: \mathbb{C}^n\to \mathbb{C}\; \mbox{holomorphic} 
%\;  
%(\forall i\in \{1, \ldots, n\}), \\ 
(*) \mbox{ is satisfied},  \\ 
Reg(\theta), Reg(b) \mbox{ are dense in }\mathbb{C}^n
\end{array}
\!
\right.\right\},     
\]
%\begin{eqnarray*}
%\mathcal{GX}_{\mathbb{C}}  & = &   
%\left\{\left.
%\Phi=\left(\left(\theta_1, \ldots, \theta_n\right), a, 
%\left(b_1, \ldots, b_n\right)\right): \mathbb{C}^n\to \mathbb{C}^{2n+1} 
%\; \mbox{ holomorphic}
%\; \right|\; 
%\theta_i, a, b_i: \mathbb{C}^n\to \mathbb{C}\quad \mbox{holomorphic}, \; 
%(*) \mbox{ in \S 1 is satisfied}\right\}, \\ 
\[
 \mathcal{GY}_{\mathbb{C}}  =  
 \left\{
\left(\left.
\left(\theta_1, \ldots, \theta_n\right),\,  
\frac{1}{2}\sum_{i=1}^n \theta_i^2, \, 
%\left(\theta_1^2+\cdots +\theta_n^2\right), 
\left(\theta_1, \ldots, \theta_n\right) 
\right) : \mathbb{C}^n\to \mathbb{C}^{2n+1} 
%\; \right|\; 
%\theta_i: \mathbb{C}^n\to \mathbb{C}
\; \mbox{ holomorphic} 
 \; \right| \; 
Reg(\theta) \mbox{ is dense in } \mathbb{C}^n 
\right\}.    
\] 
%\overline{ 
%\mbox{reg}\left(\theta_1, \ldots, \theta_n\right)} 
%= \mathbb{R}^n  
%\right\}.   
%\]
Define the \textit{complex Legendre involution} 
$\mathcal{L}_{\mathbb{C}}: 
\mathcal{GX}_{\mathbb{C}}\to \mathcal{GX}_{\mathbb{C}}$ 
by 
\[
\mathcal{L}_{\mathbb{C}}
\left(\left(\theta_1, \ldots, \theta_n\right),\; a,\;  
\left(b_1, \ldots, b_n\right)\right) 
= 
\left(\left(b_1, \ldots, b_n\right),\; \sum_{i=1}^n b_i\theta_i -a,\;  
\left(\theta_1, \ldots, \theta_n\right)\right).   
\] 
Namely, the form of the complex Legendre involution is exactly the same as 
the Legendre involution.    
Then, it is clear that the proof of Theorem \ref{theorem2} %given in this section 
simultaneously proves the following corollary.   
\begin{corollary}\label{corollary1}
Assume a transform 
$\mathcal{T}_{\C} : \mathcal{GX}_{\mathbb{C}}\to 
\mathcal{GX}_{\mathbb{C}}$ 
satisfies 
\begin{enumerate}
\item $\mathcal{T}_{\C}\circ\mathcal{T}_{\C}(\Phi) = \Phi$  
for any $\Phi\in \mathcal{GX}_{\C}$.    
\item For any $\Phi\in \mathcal{GX}_{\C}$, 
$\mathcal{T}_{\C}(\Phi)=\Phi$  if $\Phi\in \mathcal{GY}_{\mathbb{C}}$ and 
$\mathcal{T}_{\C}(\Phi) \ne \Phi$ 
if $\Phi\not\in \mathcal{GY}_{\mathbb{C}}$.     
%$\mathcal{GY}$.   
\end{enumerate}
Assume moreover that there exists a holomorphic mapping 
$F: \mathbb{C}^{2n+1}\to \mathbb{C}^{2n+1}$ satisfying 
$\mathcal{T}_{\C}(\Phi)=F\circ \Phi$ for any $\Phi\in \mathcal{GX}_{\C}$.    
Then, $\mathcal{T}_{\C}$ must be the complex Legendre involution 
$\mathcal{L}_{\mathbb{C}}: \mathcal{GX}_{\mathbb{C}}
\to \mathcal{GX}_{\mathbb{C}}$.   
\end{corollary}
%\noindent 
%Corollary \ref{corollary1} 
%suggests the existence of the complex Legendre world.     
%%%%%%%%%%%%%%%%%%%%%%%%%%%%%%%%%%%%%%%%%%%%%%%%% 
%%%%%%%%%%%%%%%%%%%%%%%%%%%%%%%%%%%%%%%%%%%%%%%%%%%%%%%%%%%%% 
%%%%%%%%%%%%%%%%%%%%%%%%%%%%%%%%%%%%%%%%%%%%%%%%%%%%%%%%%%%%%%
%\section*{Acknowledgement}
%\begin{acknowledgements}
%\end{acknowledgements}
\section*{Acknowledgements}
~The author was %partially 
supported by JSPS KAKENHI (Grant No. 23K03109).     
%\par 
%This work was supported 
%by the Research Institute for Mathematical Sciences, 
%a Joint Usage/Research Center located in Kyoto University. 
%
%
%%%%%%%%%%%%%%%%%%%%%%%%%%%%%%%%%%%%%%%%%%%%%%%%%%%%


\begin{thebibliography}{99}
\bibitem{arnoldmechanics}V.~I.~Arnol'd, 
\emph{Mathematical Methods of Classical Mechanics 2nd edition}, 
Graduate Texts in Mathematics, {\bf 60},  
Springer Netherland, Dordrecht, 1989.  
https://doi.org/10.1007/978-1-4757-2063-1 
\bibitem{arnold}V.~I.~Arnol'd, 
\emph{Singularities of Caustics and Wavefronts}, 
Mathematics and its Applications, {\bf 62}, 
Springer Netherland, Dordrecht, 1990.      
https://doi.org/10.1007/978-94-011-3330-2  
\bibitem{arnoldetal}V.~I.~Arnol'd, S.~M.~Gusein-Zade, 
and A.~N.~Varchenko, 
\emph{Singularities of Differentiable Maps I}, 
Monographs in Mathematics {\bf 82}, Birkh\"auser, 
Boston Basel Stuttgart, 1985.       
https://doi.org/10.1007/978-1-4612-3940-6    
\bibitem{artsteinavidanmilman}S.~Artstein-Avidan and 
V.~Milman, 
\emph{The concept of duality in convex analysis,
and the characterization of the Legendre transform}, 
Ann. of Math., {\bf 169} (2009), 661--674.   
https://doi.org/10.4007/annals.2009.169.661 
%\bibitem{gauss}T.~Banchoff, T.~Gaffney and C.~MacCrory, 
%\emph{Cusps of Gauss Mappings}, 
%Pitman Advanced Pub. Progman, 1982.    
%{\color{black}{
%\bibitem{bennema}P.~Bennema, 
%\emph{On the crystallographic and statistical mechanical foundations of
%the forty-year old Hartman-Perdok theory},    
%J. Crystal Growth, {\bf 166} (1996), 17--28.   
%https://doi.org/10.1016/0022-0248(96)00043-7.  
%\bibitem{blaschke}P.~Blaschke, 
%\emph{Pedal coordinates, dark Kepler and other forth problems}, 
%J.~Math.~Phys., {\bf 58} (2017), 063505.   
%https://doi.org/10.1063/1.4984905
%}}
%\bibitem{brocker}T.~Br\"{o}cker, 
%\emph{Differentiable Germs and Catastrophes}, 
%Cambridge University Press, Cambridge, 1975.   
%https://doi.org/10.1017/CBO9781107325418   
%\bibitem{brucegiblin}J.~W.~Bruce and P.~J.~Giblin,
%\emph{Curves and Singularities (second edition)}, 
%Cambridge University Press, Cambridge, 1992.      
%https://doi.org/10.1017/CBO9781139172615
%{\color{black}{
%\bibitem{brucegiblingibson}J.~W.~Bruce, P.~J.~Giblin and C.~G.~Gibson, 
%\emph{On caustics by reflection}, 
%Topology, {\bf 21} (1982), 179--199.   
%https://doi.org/10.1016/0040-9383(82)90004-0.  
%\bibitem{farouki}R.~T.~Farouki, 
%\emph{Pythagorean-Hodograph Curves: 
%Algebra and Geometry Inseparable}, 
%Geometry and Computing, {\bf 1}, Springer, 2008.   
%https://doi.org/10.1007/978-3-540-73398-0    
%\bibitem{fisher}R.~J.~Fisher and H.~Turner Laquer, 
%\emph{Hyperplane envelopes and the Clairaut equation}, 
%J.~Geom.~Anal., {\bf 20} (2010), 609--650.   
%https://doi.org/10.1007/s12220-010-9129-0    
%\bibitem{fuchs}D.~Fuchs, 
%\emph{Evolutes and involutes of spatial curves}, 
%Amer. Math. Monthly, 
%{\bf 120} (2013), 217--231
%https://doi.org/10.4169/amer.math.monthly.120.03.217   
%\bibitem{giga}Y.~Giga, 
%\emph{Surface Evolution Equations}, 
%Monographs of Mathematics, {\bf 99}, Springer, 2006.  \\  
%https://doi.org/10.1007/3-7643-7391-1
%\bibitem{feynman}D.~Goodstein and J.~R.~Goodstein, 
%\emph{Feynman's Lost Lecture:  the Motion of Planets 
%Around the Sun (first edition)}, 
%W.W. Norton \& Company, New York, 1996.   
%\bibitem{hannishimura}H.~Han and T.~Nishimura, 
%\emph{Spherical method for studying Wulff shapes and related topics}, 
%Adv. Stud. Pure Math., {\bf 78}, 1--53, Math. Soc. Japan. Tokyo, 2018.    
%https://doi.org/10.2969/aspm/07810001
%}}
%\bibitem{hoffmancahn}D.~W.~Hoffman and J.~W.~Cahn, 
%\emph{A vector thermodynamics for anisotropic surfaces}, 
%Surface Science, {\bf 31} (1972), 368--388.      
%https://doi.org/10.1016/0039-6028(72)90268-3  
\bibitem{hormander}L.~H\"{o}rmander, 
\emph{Notions of Convexity}, 
Progress in Mathematics {\bf 127},  
Birkh\"{a}user Boston, 1994.   
https://doi.org/10.1007/978-0-8176-4585-4 
%\bibitem{ishikawaaustralia}G.~Ishikawa, 
%\emph{Opening of differentiable map-germs and unfoldings}, 
%Topics on real and complex singularities, 87-–113, 
%World Sci. Publ., Hackensack, NJ, 2014.
%https://doi.org/10.1142/9789814596046\_0007
\bibitem{ishikawa}G.~Ishikawa,  
\emph{Singularities of frontals}, 
Adv. Stud. Pure Math., {\bf 78}, 
55--106, Math. Soc. Japan, Tokyo, 2018.  
https://doi.org/10.2969/aspm/07810055
%\bibitem{izumiya}S.~Izumiya, 
%\emph{Singular solutions of first order differential equations}, 
%Bull. London Math. Soc., 
%{\bf 26} (1994), 69--74.  
%https://doi.org/10.1112/blms/26.1.69
\bibitem{janeczkonishimura}S.~Janeczko and T.~Nishimura, 
\emph{ Anti-orthotomics of frontals and their applications}, J. Math. Anal. Appl., 
{\bf 487} (2020), 124019.   
https://doi.org/10.1016/j.jmaa.2020.124019 
%{\color{black}{
%\bibitem{kagatsumenishimura}D.~Kagatsume and T.~Nishimura, 
%\emph{Aperture of plane curves}, J.~Singul., 
%{\bf 12} (2015), 80--91.  
%https://doi.org/10.5427/jsing.2015.12e.  
%\bibitem{lancret}M.~A.~Lancret, 
%\emph{M\'{e}moire sur les courbes \`{a} double courbure}, 
%M\'{e}moires pr\'{e}sent\'{e}s \`{a} l’Institut,  {\bf 1} (1806), 416--454.
%\bibitem{hedgehog}R.~Langevin, G.~Levitt and H.~Rosenberg,  
%\emph{H\'erissons et multih\'erissons (enveloppes 
%parametr\'ees par leur application de Gauss)},   
%In Singularities (Warsaw, 1985), pp. 245–-253, 
%Banach Center Publ., {\bf 20}, PWN, Warsaw, 1988.  
%\bibitem{montgomery}R.~Montgomery, 
%\emph{A tour of Subriemannian Geometries, Their Geodesics and Applications}, 
%Mathematical Surveys and Monographs {\bf 91}, AMS,  
%Rhode Island Providence, 
%2002.  
%https://doi.org/http://dx.doi.org/10.1090/surv/091.    
%}}
%\bibitem{nishimura}T.~Nishimura, 
%\emph{Jacobian-squared function-germs}, Pure Appl. Math. Q., 
%{\bf 13} (2017), 711-728. 
%http://dx.doi.org/10.4310/PAMQ.2017.v13.n4.a5.   
%{\color{black}{
%\bibitem{chaos}T.~Nishimura, 
%\emph{Kato's chaos created by quadratic mappings associated 
%with spherical orthotomic curves}, J.~Singul., {\bf 21} (2020), 205--211.   
%https://doi.org/10.5427/jsing.2020.21l.   
%}}
\bibitem{nishimura}T.~Nishimura, 
\emph{Hyperplane families creating envelopes}, 
Nonlinearity, {\bf 35} (2022), 2588–-2621.   
https://doi.org/10.1088/1361-6544/ac61a0    
\bibitem{nishimurahmj}T.~Nishimura, 
\emph{Envelopes of straight line families in the plane}, 
to appear in Hokkaido Mathematical Journal 
(available at arXiv:2307.07232 [math.DG]). 
\bibitem{rockafellar}R.~Tyrrell~Rockafellar, 
\emph{Convex Analysis}, 
Princeton Landmarks in Mathematics and Physics {\bf 11}, 
Princeton University Press, 1996.    
https://doi.org/10.1515/9781400873173 
%\bibitem{sajiumeharayamada}K.~Saji, M.~Umehara and 
%K.~Yamada, 
%\emph{The geometry of fronts}, Ann. Math., 
%{\bf 169-2} (2009), 491--529.   \\ 
%https://doi.org/10.4007/annals.2009.169.491
%\bibitem{takahashi}M.~Takahashi, 
%\emph{Envelopes of families of Legendre mappings in the unit 
%tangent bundle over the Euclidean space}, 
%J. Math. Anal. Appl., {\bf 473} (2019), 408--420.        
%https://doi.org/10.1016/j.jmaa.2018.12.057. 
%{\color{black}{
%\bibitem{taylor}J.~E.~Taylor, 
%\emph{Crystalline variational problems}, Bull. Amer. Math. Soc., 
%{\bf 84} (1978), 568--588.    
%\bibitem{uribe}R.~Uribe-Vargas, 
%\emph{On Vertices, focal curvatures and
%differential geometry of space curves}, Bull. Braz. Math. Soc., New Series, 
%{\bf 36} (2005), 285--307.   
%https://doi.org/10.1007/s00574-005-0040-4.    
%\bibitem{wangnishimura}Y.~Wang and T.~Nishimura, 
%\emph{Envelopes created by circle families in the plane}, J.~Geom., 
%{\bf 115} (2024), article number 7.    
%https://doi.org/10.1007/s00022-023-00708-z 
%\bibitem{wulff}G.~Wulff, 
%\emph{Zur frage der geschwindindigkeit 
%des wachstrums und der aufl\"osung der krystallflachen}, 
%Z. Kristallographine und Mineralogie, {\bf 34} (1901), 449--530.   
%}}
\end{thebibliography}
\end{document}